\documentclass{amsart}
\usepackage{amsmath}
\usepackage{latexsym,ulem}
\usepackage{amssymb}
\usepackage{graphics}

\input epsf%
\newtheorem{theorem}{Theorem}[section]

\newtheorem{lemma}{Lemma}[section]
\newtheorem{remark}{Remark}[section]

\newtheorem{Definition}{Definition}[section]
\def\qed{{\hfill{\vrule height4pt width3pt depth2pt}}}

   \long\def\comment#1{}

\def\ad#1{\begin{aligned}#1\end{aligned}}  \def\b#1{\mathbf{#1}}
\def\a#1{\begin{align*}#1\end{align*}} \def\an#1{\begin{align}#1\end{align}}
\def\e#1{\begin{equation}#1\end{equation}} \def\d{\operatorname{div}}
\def\p#1{\begin{pmatrix}#1\end{pmatrix}} 
  \numberwithin{equation}{section}
\numberwithin{table}{section}
\numberwithin{figure}{section}

\def\boxit#1{\vbox{\hrule height1pt \hbox{\vrule width1pt\kern1pt
     #1\kern1pt\vrule width1pt}\hrule height1pt }}
 \def\lab#1{\boxit{\small #1}\label{#1}}
  \def\mref#1{\boxit{\small #1}\ref{#1}}
 \def\meqref#1{\boxit{\small #1}\eqref{#1}}

  \def\lab#1{\label{#1}} \def\mref#1{\ref{#1}} \def\meqref#1{\eqref{#1}}
\begin{document}

\newcommand{\disp}{\displaystyle}
\newcommand{\eps}{\varepsilon}
\newcommand{\To}{\longrightarrow}
\newcommand{\C}{\mathcal{C}}
\newcommand{\K}{\mathcal{K}}
\newcommand{\T}{\mathcal{T}}
\newcommand{\bq}{\begin{equation}}
\newcommand{\eq}{\end{equation}}
\long\def\comments#1{ #1}
\comments{ }

\title[Finite element  approximations]{Finite element  approximations of symmetric tensors on simplicial grids in $\mathbb{R}^n$: the lower order case
    }
\date{}

\author {Jun Hu}
\address{LMAM and School of Mathematical Sciences, Peking University,
  Beijing 100871, P. R. China.  hujun@math.pku.edu.cn}

\author {Shangyou Zhang}
\address{Department of Mathematical Sciences, University of Delaware,
    Newark, DE 19716, USA.  szhang@udel.edu }

\thanks{The first author was supported by  the NSFC Projects 11271035,  91430213 and 11421101.}

\begin{abstract}
  \vskip 15pt
In this paper, we  construct, in a unified  fashion,
lower order finite element subspaces of  spaces of symmetric
tensors with square-integrable divergence on a  domain in any dimension.
These subspaces are essentially the symmetric $H(\d)-P_{k}$ ($1\leq k\leq n$)
 tensor spaces, enriched, for each  $n-1$ dimensional simplex,
 by  $\frac{(n+1)n}{2}$ $H(\d)-P_{n+1}$  bubble functions when
 $1\leq k\leq  n-1$, and by  $\frac{(n-1)n}{2}$ $H(\d)-P_{n+1}$  bubble functions when
 $k= n$.  These spaces can be used to approximate the symmetric matrix
  field in a mixed formulation problem where the  other variable is approximated by
    discontinuous piecewise $P_{k-1}$ polynomials.   This in particular leads to first order mixed elements on  simplicial grids  with
 total degrees of freedom per element $18$ plus $3$ in 2D, 48 plus 6 in 3D.
    The previous record of  the degrees  of freedom of first order mixed elements  is,
      21 plus 3 in 2D, and   156 plus 6 in 3D, on  simplicial grids.
   We  also derive, in a unified way and without using any tools like differential forms,
  a family of auxiliary mixed finite elements in any dimension. One  example in
   this family is the Raviart-Thomas elements     in one dimension, the second
    example is   the mixed finite elements  for linear elasticity in two dimensions due to
   Arnold and  Winther,  the third   example is the mixed finite elements for linear elasticity
    in three dimensions due to Arnold, Awanou and Winther.

\noindent{\bf Keywords.}{
     mixed finite element, symmetric finite element, first order system,
     conforming finite element, simplicial grid, inf-sup condition.}

 \vskip 15pt

\noindent{\bf AMS subject classifications.}
    { 65N30, 73C02.}

\end{abstract}
\maketitle

\section{Introduction}
The constructions,  using polynomial shape functions,  of stable  pairs of finite element spaces for
approximating the pair of spaces $H(\d,\Omega; \mathbb {S})\times L^2(\Omega; \mathbb{R}^n)$
 in first order
 systems  are  a long-standing, challenging and open problem, see \cite{Arnold2002,Arnold-Awanou-Winther}.
  For mixed finite elements of linear elasticity,
     many   mathematicians have been working on this problem and  compromised
    to weakly symmetric or composite elements,  cf.
    \cite{Amara-Thomas, Arnold-Brezzi-Douglas, Arnold-Douglas-Gupta,  Johnson-Mercier,
   Morley, Stenberg-1, Stenberg-2, Stenberg-3}.  It is not until  2002 that Arnold and Winther
    were able to propose the first family of mixed finite element spaces with polynomial shape  functions in two dimensions \cite{Arnold-Winther-conforming}.  Such a two dimensional  family was extended to
    a three dimensional family of mixed elements \cite{Arnold-Awanou-Winther}, while
       the lowest order element with $k=2$ was first proposed in \cite{Adams-Cockburn}.
We refer interested
        readers to \cite{Adams-Cockburn,Arnold-Awanou,Arnold-Awanou-Winther,
    Arnold-Winther-conforming,Awanou, Chen-Wang,
     Arnold-Winther-n, Gopalakrishnan-Guzman-n,Hu-Shi,
    Man-Hu-Shi, Yi-3D, Yi, Arnold-Falk-Winther,
    Boffi-Brezzi-Fortin, Cockburn-Gopalakrishnan-Guzman,
    Gopalakrishnan-Guzman, Guzman,Hu-Man-Zhang2014,Hu-Man-Zhang2013},  for recent progress on mixed finite elements
      for linear  elasticity.

In very recent papers \cite{Hu-Zhang2014a} and \cite{Hu-Zhang2014b},  Hu and Zhang
attacked this open and challenging  problem by  initially proposing   new ideas to design
 discrete stress spaces and analyze the discrete inf-sup condition. In particular, they were able
  to construct  suitable $H(\d,\Omega;\mathbb {S})-P_k$ space, with $k\geq 3$ for 2D,
    and  $k\geq 4$ for 3D,   finite element spaces for the stress discretization in both two and three dimensions.
 In \cite{Hu2014}, Hu constructed,  in a unified fashion,  suitable $H(\d,\Omega;\mathbb {S})-P_k$ space  with
      $k\geq n+1$, and  proposed a set of degrees of freedom for the shape function space, in any dimension.

  The purpose of this paper is to extend those elements in \cite{Hu2014}
  to  lower order cases where $1\leq k\leq n$.
   Since  it is, at moment, very difficult to prove that  the pair of
    $H(\d)-P_{k}$ and $L^2-P_{k-1}$  spaces is stable,
    the $H(\d)-P_k$ space has  to be enriched by some higher order polynomials whose divergence
    are in $P_{k-1}$.
     Thanks to  \cite{Hu2014},   it suffices to control the  piecewise rigid motion space.  Hence,  we only need
       to add,   for each  $n-1$ dimensional simplex,
       $\frac{(n+1)n}{2}$ $H(\d)-P_{n+1}$ bubble functions when
     $2\leq k\leq n-1$, and   $\frac{(n-1)n}{2}$ $H(\d)-P_{n+1}$
      bubble functions when $k= n$.
      This in particular leads to first order mixed elements on  simplicial grids  with
 total degrees of freedom per element $18$ plus $3$ in 2D, 48 plus 6 in 3D.
    The previous record of  the degrees  of freedom of first order mixed elements  is,
      21 plus 3 in 2D, and   156 plus 6 in 3D, on  simplicial grids.
    These enriched bubble functions belong to the lowest order space  from a family of  auxiliary
      discrete stress spaces which, together with the  $P_{k-1}$  space, form  a stable pair of spaces for
       first order  systems. Note that these  spaces in this auxiliary family are constructed in a unified and direct way
       and that no  tools like differential forms are used. One  example in
   this auxiliary family is the Raviart--Thomas elements  in one dimension, the second
    example is   the mixed finite elements  for linear elasticity in two dimensions due to
   Arnold and  Winther \cite{Arnold-Winther-conforming},  the third   example is the mixed finite elements for linear elasticity    in three dimensions due to Arnold, Awanou and Winther \cite{Arnold-Awanou-Winther}.

We end this section by introducing  first order systems and related notations. We consider  mixed finite element methods of
 first order systems with symmetric tensors:
Find $(\sigma,u)\in\Sigma\times V :=H({\rm div},\Omega;
    \mathbb {S})
        \times L^2(\Omega; \mathbb{R}^n)$, such that
\an{\left\{ \ad{
  (A\sigma,\tau)+({\rm div}\tau,u)&= 0 && \hbox{for all \ } \tau\in\Sigma,\\
   ({\rm div}\sigma,v)&= (f,v) &\qquad& \hbox{for all \ } v\in V. }
   \right.\lab{eqn1}
}
Here the symmetric tensor space for the stress $\Sigma$ is  defined by
  \an{   \lab{S}
  H({\rm div},\Omega; \mathbb {S})
    &:= \Big\{ \tau=\p{\tau_{11} &\cdots &\tau_{1n} \\
        \vdots  & \vdots & \vdots \\
        \tau_{n1} &\cdots   & \tau_{nn}  }
     \in H(\d, \Omega; \mathbb{R}^{n\times n})
    \ \Big| \ \tau ^T = \tau  \Big\},}
and  the   space for the vector displacement  $V$ is
  \an{  \lab{V}
     L^2(\Omega; \mathbb{R}^n) &:=
     \Big\{ \p{u_1, &\cdots, &u_n}^T
          \ \Big| \ u_i \in L^2(\Omega), i=1, \cdots, n \Big\}  .}
	This paper denotes by $H^k(T; X)$ the Sobolev space consisting of
functions with domain $T\subset\mathbb{R}^n$, taking values in the
finite-dimensional vector space $X$, and with all derivatives of
order at most $k$ square-integrable. For our purposes, the range
space $X$ will be either $\mathbb{S},$ $\mathbb{R}^n,$ or
$\mathbb{R}$. Let $\|\cdot\|_{k,T}$ be the norm of $H^k(T)$, $\mathbb{S}$ denote
the space of symmetric tensors, $H({\rm div},T; \mathbb{S})$
consist of square-integrable symmetric matrix fields with
square-integrable divergence. The $H(\d)$ norm is defined by
$$\|\tau\|_{H({\rm div},T)}^2:=\|\tau\|_{0, T}^2
   +\|{\rm div}\tau\|_{0, T}^2.$$
$L^2(T; \mathbb{R}^n)$ is the space of vector-valued functions
which are square-integrable.
Here, the compliance tensor
$A=A(x):\mathbb{S} ~\rightarrow~\mathbb{S}$, characterizing the
properties of the material, is bounded and symmetric positive
definite uniformly for $x\in\Omega$.

The rest of the paper is organized as follows.  In the next section, we present
 some preliminary results from \cite{Hu2014}; see also \cite{Hu-Zhang2014a} and
 \cite{Hu-Zhang2014b}, for the cases $n=2$ and $n=3$, respectively.  In Section 3, based on these preliminary
  results, we  propose a family of  auxiliary mixed finite elements.  In Section 4, we present
   lower order mixed finite elements and  analyze the well--posedness of the discrete problem and error estimates of the approximation solution.
   In Section 5,  we present the first order mixed  elements.  In Section 6, we give a way to construct these added bubble functions  for each face in three dimensions.
 The paper ends with Section 6 which lists some numerics.

\section{Preliminary results}
Suppose that the  domain $\Omega$ is subdivided by a family of shape regular simplicial
 grids  $\mathcal{T}_h$ (with the grid size $h$).
We introduce the finite element space of order $k\geq 1$  on $\mathcal{T}_h$.
\an{\lab{Sh}
 \Sigma_{k,h} :=\Big\{ \sigma
    &\in H(\d,\Omega; \mathbb{S}),   \sigma|_K\in P_k(K;  \mathbb{S})
     \   \forall K\in\mathcal{T}_h \Big\}, }
     where $P_k(K;  \mathbb{X})$ denotes the space  of polynomials of degree $\leq k$, taking value in the space $X$.

    To define the degrees of freedom  for the shape function space $P_k(K;  \mathbb{S})$, let $\b x_0, \cdots, \b x_n$  be the  vertices of
    simplex  $K$.
The referencing mapping is then
\a{  \b x: &= F_K(\hat {\b x})
     = \b x_0 + \p{ \b x_1 -\b x_0, &
                     \cdots,&
                     \b x_n -\b x_0  } \hat{\b x}, }
      mapping the reference tetrahedron $\hat K:=\{ 0 \le \hat x_1, \cdots, \hat x_n,
          1-\sum\limits_{i=1}^n\hat x_i \le 1 \}$ to $K$.
   Then the inverse mapping is
\an{ \lab{imap}  \hat {\b x}: &= \p{\b \nu_1^T\\ \vdots \\ \b \nu_n^T} (\b x-\b x_0), }
where
 \an{\lab{nnt} \p{\b \nu_1^T\\ \vdots \\ \b \nu_n^T }=\p{ \b x_1 -\b x_0, &
      \cdots,  &
                     \b x_n -\b x_0}^{-1}. }
By \meqref{imap}, these normal vectors
        are coefficients of the barycentric variables:
  \a{ \lambda_1: &=  \b \nu_1 \cdot (\b x-\b x_0), \\
     \vdots & \\
      \lambda_n: &= \b \nu_n \cdot (\b x-\b x_0), \\
    \lambda_0: &= 1-\sum\limits_{i=1}^n\lambda_i. }
For any edge $\b x_i\b x_j$ of element $K$, $i\not =j$, let $\b t_{i,j}$ denote  associated tangent vectors, which allow for us to
 introduce the following symmetric matrices of rank one
 \begin{equation}\label{Tnew}
 T_{i,j}:=\b t_{i,j}\b t_{i,j}^T,  0\leq i< j\leq n.
 \end{equation}
For these matrices of rank one, we have the following  result from \cite{Hu2014}; see also \cite{Hu-Zhang2014a} and
 \cite{Hu-Zhang2014b}, for the cases $n=2$ and $n=3$, respectively.
 \begin{lemma}\lab{lemma2.1}
   The $\frac{(n+1)n}{2}$ symmetric tensors $T_{i,j}$ in \meqref{Tnew} are linearly independent, and
    form a basis of $\mathbb{S}$.
\end{lemma}

With these symmetric matrices $T_{i,j}$ of rank one,  we define a $H(\d,K;  {\mathbb S})$ bubble function space
\an{\lab{bh}  \Sigma_{K, k,  b}: = \sum\limits_{0\leq i<j\leq n} \lambda_i\lambda_jP_{k-2}(K;  \mathbb{R})T_{i, j}}
Define the full $H(\d,K;  {\mathbb S})$ bubble function space consisting of polynomials of degree $\leq k$
\begin{equation}
\Sigma_{\partial K,  k, 0}:=\{\tau\in H(\d, K;  \mathbb{S})\cap P_{k}(K; \mathbb{S}), \tau\b \nu|_{\partial K}=0 \}.
\end{equation}
Here $\nu$ is the   normal vector of $\partial K$.

 We have the following  result  due to \cite{Hu2014}.

   \begin{lemma}\label{lemma2.2}  It holds that
   \begin{equation}
   \Sigma_{K, k, b}=\Sigma_{\partial K,  k, 0}.
   \end{equation}
      \end{lemma}

We need  an important result concerning the divergence space of the bubble function space.
 To this end, we introduce the following rigid motion space on each element $K$.
  \begin{equation}\label{R}
  R(K):=\{ v\in H^1(K; \mathbb{R}^n), (\nabla v+\nabla v^T)/2=0\}.
  \end{equation}
It follows from the definition that  $R(K)$ is a subspace of $P_1(K;  \mathbb{R}^n)$.
For $n=1$, $R(K)$ is the constant function space over $K$. The dimension of $R(K)$ is $\frac{n(n+1)}{2}$.
For two dimensions, the rigid motion space $R(K)$ is
\begin{equation}\label{2Drigid}
R(K):=\bigg\{\begin{pmatrix} a_1 \\a_2 \end{pmatrix}+ b\begin{pmatrix}-x_2\\ x_1\end{pmatrix}, a_1, a_2, b\in \mathbb{R}\bigg\};
\end{equation}
for three dimensions, the rigid motion space $R(K)$ reads
\begin{equation}\label{3Drigid}
R(K): =\bigg\{\begin{pmatrix} a_1 \\a_2  \\ a_3\end{pmatrix}+ b_1\begin{pmatrix}-x_2\\ x_1\\0\end{pmatrix}+b_2\begin{pmatrix}-x_3\\ 0\\x_1\end{pmatrix}+b_3\begin{pmatrix}0\\ -x_3\\x_2\end{pmatrix}, a_i, b_i\in \mathbb{R},  i=1,2,3 \bigg\}.
\end{equation}
This allows for defining the orthogonal complement  space of  $R(K)$ with  respect to $P_{k-1}(K; \mathbb{R}^n)$ by
 \begin{equation}
 R^\perp(K):=\{v\in P_{k-1}(K; \mathbb{R}^n), (v, w)_K=0\text{ for any }w\in R(K)\},
 \end{equation}
 where the inner product $(v, w)_K$ over $K$ reads $(v, w)_K=\int_K v\cdot w d\b x$.
\begin{lemma}\label{lemma2.3} For any $K\in\mathcal{T}_h$, it holds that
\begin{equation}
\d \Sigma_{K, k, b}=R^\perp(K).
\end{equation}
\end{lemma}
\begin{proof} The proof can be found in \cite{Hu2014}; see also \cite{Hu-Zhang2014a} and
 \cite{Hu-Zhang2014b}, for the cases $n=2$ and $n=3$, respectively.
\end{proof}

We need a classical result and its variant.
\begin{lemma}\label{lemma2.4} It holds the following  Chu-Vandermonde combinatorial identity and its variant
\begin{equation}\label{van}
\sum\limits_{\ell=0}^nC_{n+1}^{\ell+1}C_{k-1}^\ell =\sum\limits_{\ell=0}^nC_{n+1}^{n-\ell}C_{k-1}^\ell=C_{n+k}^n,
\end{equation}
and
\begin{equation}\label{varvan}
\sum\limits_{\ell=0}^nC_{n+1}^{\ell+1}C_{k-1}^\ell C_{\ell+1}^2=\frac{(n+1)n}{2}C_{n+k-2}^n,
\end{equation}
where the combinatorial number $C_{n}^m=\frac{n\cdots (n-m+1)}{m\cdots 1}$ for $n\geq m$ and $C_{n}^m=0$ for $n<m$.
\end{lemma}

\section{A family of auxiliary mixed  elements in any dimension}
\subsection{The lowest order auxiliary mixed  elements}
To define  lower order  mixed finite elements with $k\leq n$, we need a family of auxiliary mixed  elements.  To this end,
 we  introduce the following divergence free space for element $K\in \mathcal{T}_h$,
 \begin{equation}
 \Sigma_{3\rightarrow n+1, DF}(K; \mathbb{S}):=\{\tau\in P_{n+1}(K; \mathbb{S})\backslash P_2(K; \mathbb{S}), \d \tau=0 \}.
 \end{equation}
 It is straightforward to see that the dimension of the space $\Sigma_{3\rightarrow n+1, DF}(K; \mathbb{S})$ reads
 \begin{equation}
\bigg( \frac{(2n+1)!}{n!(n+1)!}-\frac{(n+2)!}{2!n!}\bigg)\frac{n(n+1)}{2}-n\frac{(2n)!}{n!n!}+n(n+1).
 \end{equation}
 Here $\bigg( \frac{(2n+1)!}{n!(n+1)!}-\frac{(n+2)!}{2!n!}\bigg)\frac{n(n+1)}{2}$ is the dimension of the space $P_{n+1}(K,\mathbb{S})\backslash P_2(K, \mathbb{S})$, and $n\bigg(\frac{(2n)!}{n!n!}-(n+1)\bigg)$ is the number of constraints
  by the divergence free. Then we can define the following enriched $P_2(K;  \mathbb{S})$ space
  \begin{equation}
  P_2^\ast(K; \mathbb{S}):=P_2(K; \mathbb{S})+\Sigma_{3\rightarrow n+1, DF}(K; \mathbb{S}).
  \end{equation}
   It follows that the dimension of $P_2^\ast(K; \mathbb{S})$ is equal to
   \begin{equation}\label{eq3.4}
   \begin{split}
    \text{the~~ dimension~~ of ~~~}P_2(K; \mathbb{S})+ \text{~~~ the~~~ dimension~~~ of~~ }\Sigma_{3\rightarrow n+1, DF}(K;  \mathbb{S})\\
    = \frac{(2n+1)!}{n!(n+1)!}\frac{n(n+1)}{2}-n\frac{(2n)!}{n!n!}+n(n+1).
   \end{split}
   \end{equation}
  To  present the degrees of freedom of  $P_2^\ast(K; \mathbb{S})$, we define
  \begin{equation}
    M_2(K):=\{\tau\in P_2^\ast(K;  \mathbb{S}), \d\tau=0 \text{ and }\tau\nu|_{\partial K}=0 \},
  \end{equation}
  where $\nu$ is the   normal vector of $\partial K$. For the space $M_2(K)$, we have the following important result.
  \begin{lemma} The dimension of $M_2(K)$ is
  \begin{equation}
  \frac{(2n-1)!}{n!(n-1)!}\frac{n(n+1)}{2}+\frac{n(n+1)}{2}-n\frac{(2n)!}{n!n!}.
  \end{equation}
  \end{lemma}
  \begin{proof} The dimension of the space $\Sigma_{K, n+1, b}$ reads
  \begin{equation}
  \frac{(2n-1)!}{n!(n-1)!}\frac{n(n+1)}{2}.
  \end{equation}
  Since the dimension of $R(K)$ is $\frac{n(n+1)}{2}$, the dimension of $R^\perp(K)$ (with respect to $P_{n}(K, \mathbb{R}^n)$) is
  \begin{equation}
  n\frac{(2n)!}{n!n!}-\frac{n(n+1)}{2}.
  \end{equation}
  It follows from the definition of $P_2^\ast(K; \mathbb{S})$ and Lemma \ref{lemma2.2} that $M_2(K)$ contains all divergence free  tensor-value functions of  $\Sigma_{K, n+1, b}$. Then the desired result follows from Lemma \ref{lemma2.3}.
  \end{proof}

\begin{theorem}\label{theorem3.1} A  matrix field  $\tau\in P_2^\ast(K; \mathbb{S})$ can be uniquely determined by the following degrees of
 freedom:
\begin{enumerate}
 \item For each $\ell$ dimensional simplex $\triangle_\ell$ of $K$, $0\leq \ell\leq n-1$,  with $\ell$ linearly independent tangential vectors $\b t_1, \cdots, \b t_\ell$, and $n-\ell$ linearly  independent normal vectors $\b \nu_1, \cdots, \b \nu_{n-\ell}$, the mean moments of degree at most $n-\ell$ over $\triangle_\ell$, of~~ $\b t_l^T\tau \b \nu_i$, $\b \nu_i^T\tau\b \nu_j$, $l=1, \cdots, \ell$, $i,j=1, \cdots, n-\ell$,  $\big( C_{n+1-\ell}^2+\ell (n-\ell)\big) C_{n}^\ell=\frac{(n-\ell)(n+\ell+1)}{2}C_{n}^\ell$ degrees of freedom  for each $\triangle_\ell$;
 \item the average of $\tau$ over $K$,  $\frac{n(n+1)}{2}$  degrees of freedom;
 \item the values of moments $\int_K \tau:\theta d\b x $, $\theta \in M_2(K)$, $\frac{(2n-1)!}{n!(n-1)!}\frac{n(n+1)}{2}+\frac{n(n+1)}{2}-n\frac{(2n)!}{n!n!}$ degrees of freedom.
 \end{enumerate}
  \end{theorem}
\begin{proof} We assume that all degrees of freedom vanish and show that $\tau=0$. Note that the mean moment becomes the value of $\tau$ for a $0$ dimensional simplex $\triangle_0$, namely, a vertex, of $K$.  The  first set of degrees of freedom implies that $\tau \b \nu =0$ on $\partial K$ while the second set of degrees of freedom shows $\d \tau=0$. Then the third set of
 degrees of freedom  proves that $\tau=0$. Next we shall
 prove that the sum of these  degrees of freedom is equal to the dimension of the space $P_2^\ast(K, \mathbb{S})$. In fact
  the sum of the first set of degrees of freedom is
  \begin{equation*}
  \sum\limits_{\ell=0}^{n-1}C_{n+1}^{\ell+1}\frac{(n-\ell)(n+\ell+1)}{2}C_{n}^\ell,
  \end{equation*}
  we refer interested readers to \cite[Theorem 2.1]{Hu2014} for a detailed proof of the numbers of degrees of freedom in the first set.
  By the Chu-Vandermonde combinatorial identity \eqref{van} and its variant \eqref{varvan}, see more details from \cite{Hu2014},
  \begin{equation*}
  \sum\limits_{\ell=0}^{n-1}C_{n+1}^{\ell+1}\frac{(n-\ell)(n+\ell+1)}{2}C_{n}^\ell=\frac{(2n+1)!}{n!(n+1)!}\frac{n(n+1)}{2}-
  \frac{(2n-1)!}{n!(n-1)!}\frac{n(n+1)}{2}.
  \end{equation*}
  Hence the  desired result follows from \eqref{eq3.4}, and the sum of the second and third sets of degrees of freedom.
\end{proof}

Then we define
\begin{equation}\label{Shplus}
\Sigma_{2,h}^\ast:=\{\tau\in H(\d, \Omega; \mathbb{S}), \tau|_K\in P_2^\ast(K; \mathbb{S}) \text{ for  any }K\in\mathcal{T}_h\}.
\end{equation}

\begin{remark} For $n=2$, we  recover the  lowest order element in \cite{Arnold-Winther-conforming};
 for $n=3$ we obtain the lowest order element
  in \cite{Adams-Cockburn}, see also \cite{Arnold-Awanou-Winther}.
\end{remark}

To define a  family of first order mixed elements, we need a   family of  simplified lowest order mixed elements, which is defined by
 \begin{equation}
  \hat{P}_2^\ast(K; \mathbb{S}):=\{\tau\in P_2^\ast(K; \mathbb{S}), \d\tau\in R(K)\}.
  \end{equation}
  The dimension of  $\hat{P}_2^\ast(K; \mathbb{S})$ is
  $$
  \frac{(2n+1)!}{n!(n+1)!}\frac{n(n+1)}{2}-n\frac{(2n)!}{n!n!}+\frac{n(n+1)}{2}.
  $$
 A complete set of degrees of freedom for $\hat{P}_2^\ast(K; \mathbb{S})$  is obtained by removing the $\frac{n(n+1)}{2}$ average values over $K$
  for   $P_2^\ast(K; \mathbb{S})$.
Then we define
\begin{equation}\label{hatShplus}
\hat{\Sigma}_{2,h}^\ast:=\{\tau\in H(\d, \Omega; \mathbb{S}), \tau|_K\in \hat{P}_2^\ast(K; \mathbb{S}) \text{ for  any }K\in\mathcal{T}_h\}.
\end{equation}
\begin{remark} For $n=2, 3$, we  recover the simplified   lowest order elements in \cite{Arnold-Winther-conforming} and  \cite{Arnold-Awanou-Winther}, respectively.
\end{remark}

\subsection{Higher order auxiliary mixed  elements}
To define auxiliary mixed elements of order $k>2$,  we  introduce the following divergence free space for element $K\in \mathcal{T}_h$,
 \begin{equation}
 \Sigma_{k+1\rightarrow k+n-1, DF}(K; \mathbb{S}):=\{\tau\in P_{k+n-1}(K; \mathbb{S})\backslash P_k(K; \mathbb{S}), \d \tau=0 \}.
 \end{equation}
Since the dimension of the space $P_{k+n-2}(K; \mathbb{R})\backslash P_{k-1}(K; \mathbb{R})$ is
$$
\frac{((k+2n-2))!}{n!(k+n-2)!}-\frac{(n+k-1)!}{n!(k-1)!},
$$
the number of the divergence free constraints is
$$
n\bigg(\frac{((k+2n-2))!}{n!(k+n-2)!}-\frac{(n+k-1)!}{n!(k-1)!}\bigg).
$$
In addition,  the dimension of the space $P_{k+n-1}(K; \mathbb{S})\backslash P_k(K; \mathbb{S})$ is
$$
\bigg( \frac{(k+2n-1)!}{n!(k+n-1)!}-\frac{(n+k)!}{k!n!}\bigg)\frac{n(n+1)}{2}.
$$
It follows that the dimension of the space $\Sigma_{k+1\rightarrow k+n-1, DF}(K; \mathbb{S})$ is
 \begin{equation}
\bigg( \frac{(k+2n-1)!}{n!(k+n-1)!}-\frac{(n+k)!}{k!n!}\bigg)\frac{n(n+1)}{2}-n\bigg(\frac{((k+2n-2))!}{n!(k+n-2)!}-\frac{(n+k-1)!}{n!(k-1)!}\bigg).
 \end{equation}
Define the following enriched $P_k(K; \mathbb{S})$  space
  \begin{equation}
  P_k^\ast(K; \mathbb{S}):=P_k(K; \mathbb{S})+\Sigma_{k+1\rightarrow k+n-1, DF}(K; \mathbb{S}).
  \end{equation}
   It follows that the dimension of $P_k^\ast(K; \mathbb{S})$ is equal to
   \begin{equation}\label{eq3.15}
   \begin{split}
    \text{the~~ dimension~~ of ~~~}P_k(K; \mathbb{S})+ \text{~~~ the~~~ dimension~~~ of~~ }\Sigma_{k+1\rightarrow k+n-1, DF}(K; \mathbb{S})\\
    =  \frac{(k+2n-1)!}{n!(k+n-1)!}\frac{n(n+1)}{2}-n\bigg(\frac{(k+2n-2)!}{n!(k+n-2)!}-\frac{(n+k-1)!}{n!(k-1)!}\bigg).
   \end{split}
   \end{equation}
  To  present the degrees of freedom of  $P_k^\ast(K; \mathbb{S})$, we define
  \begin{equation}
    M_k(K):=\{\tau\in P_k^\ast(K;  \mathbb{S}), \d\tau=0 \text{ and }\tau\nu|_{\partial K}=0 \},
  \end{equation}
  where $\nu$ is the   normal vector of $\partial K$. For the space $M_k(K)$, we have the following important result.
  \begin{lemma} The dimension of $M_k(K)$ is
  \begin{equation}\label{eq3.17}
   \frac{(k+2n-3)!}{n!(k+n-3)!}\frac{n(n+1)}{2}+\frac{n(n+1)}{2}-n\frac{(k+2n-2)!}{n!(k+n-2)!}.
  \end{equation}
  \end{lemma}
  \begin{proof} The dimension of the space $\Sigma_{K, k+n-1, b}$ reads
  \begin{equation}
  \frac{(k+2n-3)!}{n!(k+n-3)!}\frac{n(n+1)}{2}.
  \end{equation}
  Since the dimension of $R(K)$ is $\frac{n(n+1)}{2}$, the dimension of $R^\perp(K)$ (with respect to $P_{k+n-2}(K;  \mathbb{R}^n)$) is
  \begin{equation}
  n\frac{(k+2n-2)!}{n!(k+n-2)!}-\frac{n(n+1)}{2}.
  \end{equation}
  It follows from the definition of $P_k^\ast(K; \mathbb{S})$ and Lemma \ref{lemma2.2} that $M_k(K)$ contains all divergence free  tensor-value functions of  $\Sigma_{K, k+n-1, b}$. Then the desired result follows from Lemma \ref{lemma2.3}.
  \end{proof}
\begin{theorem} A  matrix field  $\tau\in P_k^\ast(K; \mathbb{S})$ can be uniquely determined by the following degrees of
 freedom:
\begin{enumerate}
 \item For each $\ell$ dimensional simplex $\triangle_\ell$ of $K$, $0\leq \ell\leq n-1$,  with $\ell$ linearly independent tangential vectors $\b t_1, \cdots, \b t_\ell$, and $n-\ell$ linearly  independent normal vectors $\b \nu_1, \cdots, \b \nu_{n-\ell}$, the mean moments of degree at most $k+n-\ell-2$ over $\triangle_\ell$, of~~ $\b t_l^T\tau \b \nu_i$, $\b \nu_i^T\tau\b \nu_j$, $l=1, \cdots, \ell$, $i,j=1, \cdots, n-\ell$,  $\big( C_{n+1-\ell}^2+\ell (n-\ell)\big) C_{k+n-2}^\ell=\frac{(n-\ell)(n+\ell+1)}{2}C_{k+n-2}^\ell$ degrees of freedom  for each $\triangle_\ell$;
 \item the values $\int_K \tau: \theta d\b x$ for any $\theta \in \epsilon(P_{k-1}(K;  \mathbb{R}^n))$, $nC_{n+k-1}^n-\frac{n(n+1)}{2}$ degrees of freedom;
 \item the values $\int_K \tau: \theta d\b x$ for any $\theta \in M_k(K)$,  $\frac{(k+2n-3)!}{n!(k+n-3)!}\frac{n(n+1)}{2}+\frac{n(n+1)}{2}-n\frac{(k+2n-2)!}{n!(k+n-2)!}$ degrees of freedom
 \end{enumerate}
  \end{theorem}
\begin{proof} We assume that all degrees of freedom vanish and show that $\tau=0$. Note that the mean moment become the value of $\tau$ for a $0$ dimensional simplex $\triangle_0$, namely, a vertex, of $K$.  The  first set of degrees of freedom implies that $\tau \b \nu =0$ on $\partial K$ while the second set of degrees of freedom shows $\d \tau=0$. Then the third set of degrees of freedom  proves that $\tau=0$.

Next we shall
 prove that the sum of these  degrees of freedom is equal to the dimension of the space $P_k^\ast(K; \mathbb{S})$. In fact, it follows from
 the Chu-Vandermonde combinatorial identity \eqref{van} and its variant \eqref{varvan} that the number of degrees in the first
  set is
  \begin{equation}
  \sum\limits_{\ell=0}^{n-1}C_{n+1}^{\ell+1}\frac{(n-\ell)(n+\ell+1)}{2}C_{n+k-2}^\ell=\frac{n(n+1)}{2}(C_{k+2n-1}^n-C_{k+2n-3}^n),
  \end{equation}
    we refer interested readers to \cite[Theorem 2.1]{Hu2014} for a detailed proof of the numbers of degrees of freedom in the first set.
  The desired result follows from \eqref{eq3.15} and \eqref{eq3.17}.
\end{proof}
Then we define
\begin{equation}\label{Shplushorder}
\Sigma_{k,h}^\ast:=\{\tau\in H(\d, \Omega; \mathbb{S}), \tau|_K\in P_k^\ast(K; \mathbb{S}) \text{ for  any }K\in\mathcal{T}_h\}.
\end{equation}

\begin{remark} For $n=2$, we  recover the higher order elements in \cite{Arnold-Winther-conforming}; for $n=3$ we obtain the higher order elements in    \cite{Arnold-Awanou-Winther}.
\end{remark}

\section{A  family of lower order mixed elements}\label{sec3}
\subsection{Mixed methods}
We propose to use the spaces $\Sigma_{k,h}$, with $2\leq k\leq n$,  defined in \eqref{Sh} to approximate $\Sigma$.  In order to  get a stable pair of spaces, we take the discrete displacement space  as the full $C^{-1}$-$P_{k-1}$ space
 \an{ \lab{Vh}
   V_{k,h}: = \{v\in L^2(\Omega; \mathbb{R}^n),  \
         v|_K\in P_{k-1}(K;  \mathbb{R}^n)\ \hbox{ for all } K\in\mathcal{T}_h \}.
    }
Unfortunately, we can not  establish the stability of the  pair of spaces $\Sigma_{k,h}$  and $V_{k,h}$.
We propose to enrich $\Sigma_{k,h}$ by some $n-1$ dimensional simplex bubble  function spaces.
Given a $n-1$ dimensional simplex $F$ of $\mathcal{T}_h$, let $\omega_F:=K^-\cup K^+$ denote the union of two elements that share $F$.  Define
\begin{equation}
\begin{split}
 \mathbb{B}_F^1:=&\bigg\{\tau\in \Sigma_{2,h}^\ast, \tau=0 \text{ on }\Omega\backslash\omega_F, \int_F \tau\nu\cdot p ds=0 \text{ for any } p\in \big(R(\omega_F)|_F\big)^\perp,\\
 & \text{ the averages of }\tau \text{ over both } K^- \text{ and } K^+ \text{ vanish },\\
 & \text{ the values of } \int_K\tau:\theta d\b x \text{ vanish for any }\theta\in M_2(K), K=K^- \text{ and } K^+ \bigg\}.
\end{split}
\end{equation}
Here $\nu$ is the   normal vector of $F$, and $\big(R(\omega_F)|_F\big)^\perp$ is the orthogonal complement of
the restriction  $R(\omega_F)|_F$ on $F$ of $R(\omega_F)$ with respect to the $L^2$ inner product  over $F$.
 We also need a subspace of $\mathbb{B}_F^1$ defined by
 \begin{equation}
 \mathbb{B}_F^2:=\{\tau\in \mathbb{B}_F^1, \int_F\tau\nu\cdot p d s=0 \text{ for any }p\in P_0(F, \mathbb{R}^n)\}.
 \end{equation}

 Hence we define the following enriched stress space
 \begin{equation}
   \Sigma_{k, h}^+=\Sigma_{k, h}+\sum\limits_{F} \mathbb{B}_F^1 \text{ for }2\leq k\leq n-1;
 \end{equation}
 and
 \begin{equation}
  \lab{Sh-plus}
   \Sigma_{k, h}^+=\Sigma_{k, h}+\sum\limits_{F} \mathbb{B}_F^2 \text{ for }k=n.
 \end{equation}

\begin{lemma}\label{lemma4.1} The space $\Sigma_{k, h}^+$ is a direct sum of the spaces $\Sigma_{k, h}$ and $\sum\limits_{F} \mathbb{B}_F^1$ for
$2\leq k\leq n-1$; is a direct sum of the spaces $\Sigma_{k, h}$ and $\sum\limits_{F} \mathbb{B}_F^2$ for
$k=n$.
\end{lemma}
\begin{proof} We first  prove the first part of the theorem. It suffices to show that, given $K\in \mathcal{T}_h$, assume    the following degrees of freedom vanish  for $\tau\in P_k(K, \mathbb{S})$ with $2\leq k\leq n-1$, then $\tau\nu=0$ on
  $\partial K$ where $\nu$ is the   normal vector of $\partial K$.
\begin{itemize}
 \item For each $\ell$ dimensional simplex $\triangle_\ell$ of $K$, $0\leq \ell\leq n-2$,  with $\ell$ linearly independent tangential vectors $\b t_1, \cdots, \b t_\ell$, and $n-\ell$ linearly  independent normal vectors $\b \nu_1, \cdots, \b \nu_{n-\ell}$, the mean moments of degree at most $n-\ell$ over $\triangle_\ell$, of~~ $\b t_l^T\tau \b \nu_i$, $\b \nu_i^T\tau\b \nu_j$, $l=1, \cdots, \ell$, $i,j=1, \cdots, n-\ell$,  $\big( C_{n+1-\ell}^2+\ell (n-\ell)\big) C_{n}^\ell=\frac{(n-\ell)(n+\ell+1)}{2}C_{n}^\ell$ degrees of freedom  for each $\triangle_\ell$;
 \end{itemize}
 In fact, it follows from \cite[Theorem 2.1]{Hu2014} that such a set of degrees of freedom indicates the $\tau\nu=0$ on
  $\partial K$.

  Next we turn to the second part of the theorem.  For this case, if the above set of degrees of freedom
  and  the following set of degrees of  freedom
   \begin{itemize}
   \item  the  average moment of  degree zero of $\tau \nu$ for any $n-1$ dimensional simplex $\triangle_{n-1}$ of $K$ with the   normal vector
   \end{itemize}
    vanish,  we have $\tau\nu=0$ on $\partial K$, see \cite[Theorem 2.1]{Hu2014} for more details.  This completes the proof.
\end{proof}

 It follows from the definition of $V_{k,h}$ ( $P_{k-1}$ polynomials)
    and $\Sigma_{k,h}^+$ (enriched $P_k$ polynomials) that
   \a{ \d  \Sigma_{k,h}^+ \subset V_{k,h}.}
This, in turn, leads to a strong divergence-free space:
 \an{ \lab {kernel}
    Z_h&:= \{\tau_h\in\Sigma_{k,h}^+ \ | \ (\d\tau_h, v)=0 \quad
	\hbox{for all } v\in V_{k,h}\}\\
    \nonumber
          &= \{\tau_h \in\Sigma_{k,h}^+ \ | \  \d \tau_h=0
    \hbox{\ pointwise } \}.
    }

The mixed finite element approximation of Problem (1.1) reads: Find
   $(\sigma_h,~u_h)\in\Sigma_{k,h}^+\times V_{k,h}$ such that
 \e{ \left\{ \ad{
    (A\sigma_h, \tau)+({\rm div}\tau, u_h)&= 0 &&
              \hbox{for all \ } \tau \in\Sigma_{k,h}^+,\\
     (\d\sigma_h, v)& = (f, v) &&  \hbox{for all \ } v\in V_{k,h}.
      } \right. \lab{DP}
    }

\subsection{Stability analysis and error estimates}
The convergence of the finite element solution follows
   the stability and the standard approximation property.
So we consider first the well-posedness  of the discrete problem
    \meqref{DP}.
By the standard theory,  we only need to prove
   the following two conditions, based on their counterpart at
    the continuous level.

\begin{enumerate}
\item K-ellipticity. There exists a constant $C>0$, independent of the
   meshsize $h$ such that
    \an{ \lab{below} (A\tau, \tau)\geq C\|\tau\|_{H(\d)}^2\quad
       \hbox{for all } \tau \in Z_h, }
    where $Z_h$ is the divergence-free space defined in \meqref{kernel}.

\item  Discrete B-B condition.
    There exists a positive constant $C>0$
            independent of the meshsize $h$, such that
    \an{\lab{inf-sup}
   \inf_{0\neq v\in V_{k,h}}   \sup_{0\neq\tau\in\Sigma_{k,h}^+}\frac{({\rm
        div}\tau, v)}{\|\tau\|_{H(\d)}  \|v\|_{0} }\geq
    C .}
\end{enumerate}

It follows from $\d  \Sigma_{k,h}^+ \subset V_{k,h}$ that $\d  \tau=0$ for
   any $\tau\in Z_h$. This implies the above K-ellipticity condition
	\meqref{below}.
It remains to show the discrete B-B condition \meqref{inf-sup},
  in the following two lemmas.

  For the analysis, we need a subspace $\widetilde \Sigma_{k,h}:=\Sigma_{k,h}\cap H^1(\Omega, \mathbb{S})$ of $\Sigma_{k,h}$. For $\tau\in \widetilde \Sigma_{k,h}$, the degrees of freedom on any element $K$
   are:  for each $\ell$ dimensional simplex $\triangle_\ell$ of $K$, $0\leq \ell\leq n$,  the mean moments of degree at most $k-\ell-1$ over $\triangle_\ell$, of $\tau$.  A standard argument is able to prove that these degrees of freedom are unisolvent.

\begin{lemma}\label{lemma1}
For any $v_h\in V_{k,h}$,  there is a $\tau_h \in
    \widetilde \Sigma_{k,h}+\sum\limits_{F}\mathbb{B}_F^1$  with $2\leq k\leq n-1$ such that,
  for all polynomial $p\in R(K)$, $K\in\mathcal{T}_h$,
   \bq\label{l-1}  \int_K (\d\tau_h-v_h) \cdot p\, d\b x=0
      \quad \hbox{\rm
      and } \quad \|\tau_h\|_{H(\d)}\leq C\|v_h\|_{0}. \eq
     \end{lemma}
\begin{proof}
  Let $v_h\in V_{k,h}$.
  By the stability of the continuous formulation,
      cf. \cite{Arnold-Winther-conforming} for two dimensional case, there is
    a $\tau \in  H^1(\Omega; \mathbb{S})$ such that,
   \a{ \d\tau=v_h \quad \hbox{\rm
      and } \quad \|\tau\|_{1}\le  C\|v_h\|_{0}. }
In this paper, we only consider the domain such that the above stability holds. We refer interested authors
    to \cite{GiraultRaviart1986} for the  classical result which states  it  is true for Lipschitz domains in
     $\mathbb{R}^n$; see \cite{Duran2001} for more refined results.
First let $I_h$ be a  Scott-Zhang \cite{Scott-Zhang}
   interpolation operator such that
  \begin{equation}\label{eq3.9}
  \|\tau-I_h\tau\|_0+h\|\nabla I_h\tau\|_0\leq Ch\|\nabla\tau\|_0.
  \end{equation}
   These enriched bubble  functions in $\sum\limits_{F}\mathbb{B}_F^1$  on the $n-1$ dimensional simplices $F$ allow for defining a correction $\delta_h\in \sum\limits_{F}\mathbb{B}_F^1$ such that
    \begin{equation}
    \int_{F} \delta_h\nu \cdot p d \b s=\int_{F} (\tau-I_h\tau)\nu\cdot p d\b s\text{ for any }
    p\in R(K)|_{F}.
    \end{equation}
    Finally we take
    \begin{equation}
    \tau_h=I_h\tau+\delta_h.
    \end{equation}
    We get  a partial-divergence matching property of $\tau_h$:
     for any $p\in R(K)$,  as the symmetric gradient $\epsilon( p)=0$,
 \a{  \int_K (\d\tau_h-v_h) \cdot p \, d\b x
     & = \int_K (\d\tau_h-\d \tau) \cdot p \, d\b x\\
        & = \int_{\partial K}  (\tau_h-\tau)\nu\cdot p\,  d\b s =0.}
 The  stability estimate follows from \eqref{eq3.9} and  the definition of the correction $\delta_h$.
\end{proof}

\begin{remark}\label{remark4.1} A  modification of the above proof applies for the case where $k=n$.  In fact,
the bubble functions in the spaces $\widetilde \Sigma_{k,h}$ and  $\sum\limits_{F}\mathbb{B}_F^2$
on the $n-1$ dimensional simplices $F$ are able   to control the  constant  subspace of $R(K)$
and its  orthogonal complement, respectively.
\end{remark}

 We are in the position to show the well-posedness of the discrete problem.
\begin{theorem}
 For the discrete problem (\ref{DP}), the K-ellipticity \meqref{below}
    and the discrete B-B
 condition \meqref{inf-sup} hold uniformly.
  Consequently,  the discrete
     mixed problem \meqref{DP} has a unique solution
         $(\sigma_h,~u_h)\in\Sigma_{k,h}^+\times V_{k,h}$.
\end{theorem}
\begin{proof}  The  K-ellipticity immediately follows from the fact
      that $\d  \Sigma_{k,h}^+ \subset V_{k,h}$.
      To prove the  discrete B-B  condition \meqref{inf-sup},
      for any $v_h\in V_{k,h}$,
        it follows from Lemma \ref{lemma1} and Remark \ref{remark4.1}  that there exists a
      $\tau_{1}\in \Sigma_{k,h}^+$ such that,  for any polynomial $p\in R(K)$,
   \bq  \int_K (\d\tau_1-v_h) \cdot pd\b x=0
      \quad \hbox{\rm
      and } \quad \|\tau_1\|_{H(\d)}\leq C\|v_h\|_{0}. \eq

Then it follows from  Lemma  \ref{lemma2.3} that
 there is a $\tau_2\in\Sigma_{k, h}$ such that  $\tau_2|_K \in \Sigma_{K, k, b}$  and
   \begin{equation}
     \d\tau_2 = v_h-\d\tau_1, \|\tau_2\|_0=\min\{\|\tau\|_0, \d\tau=v_h-\d\tau_1, \tau\in \Sigma_{K, k, b} \}
     \end{equation}
     It follows from the definition of $\tau_2$ that $\|\d \tau_2\|_0$  defines a norm for it. Then, a scaling argument proves
\begin{equation}
 \|\tau_2 \|_{H(\d)}\leq C\|\d\tau_1-v_h\|_{0}.
\end{equation}
Let $\tau=\tau_1+\tau_2$.  This implies that
\begin{equation}
\d\tau=v_h \text{ and } \|\tau\|_{H(\d)}\leq C\|v_h\|_{0},
\end{equation}
this proves the discrete B-B condition \meqref{inf-sup}.
\end{proof}

\begin{theorem}\label{MainError} Let
  $(\sigma, u)\in\Sigma\times V$ be the exact solution of
   problem \meqref{eqn1} and $(\tau_h, u_h)\in\Sigma_{k,h}^+\times
   V_{k,h}$ the finite element solution of \meqref{DP}.  Then,
   for $2\leq k\leq n$,
\an{ \lab{t1} \|\sigma-\sigma_h\|_{H({\rm div})}
    + \|u-u_h\|_{0}&\le     Ch^k(\|\sigma\|_{k+1}+\|u\|_{k}).
      }
\end{theorem}

\begin{proof}
 The stability of the elements and the standard theory of mixed
  finite element methods \cite{Brezzi, Brezzi-Fortin} give the
  following quasioptimal error estimate immediately
\an{
  \label{theorem-err1} \|\sigma-\sigma_h\|_{H({\rm
  div})}+\|u-u_h\|_{0}\leq C \inf\limits_{\tau_h\in\Sigma_{k,h}^+,v_h\in
  V_{k,h}}\left(\|\sigma-\tau_h\|_{H({\rm div})}+\|u-v_h\|_{0}\right).}
Let $P_h$ denote the local $L^2$ projection operator,
   or  element-wise interpolation operator,  from $V$ to $V_{k,h}$,
  satisfying the error estimate
\an{\label{proj-error}
   \|v-P_hv\|_{0}\leq Ch^k\|v\|_{k} \text{ for any }v\in H^k(\Omega; \mathbb{R}^n). }
Choosing $\tau_h=I_h\sigma\in \Sigma_{k,h}$
    where $I_h$ is defined in \eqref{eq3.9} as $I_h$ preserves symmetric $P_k$ functions locally,
   \an{ \lab{p-err2}
       \|\sigma -\tau_h\|_{0} + h |\sigma -\tau_h|_{H(\d)}
        \le Ch^{k+1} \|\sigma\|_{k+1}. }
 Let $v_h= P_h v$ and $\tau_h=I_h\sigma$ in (\ref{theorem-err1}),
  by (\ref{proj-error}) and \meqref{p-err2}, we
   obtain  \meqref{t1}.
\end{proof}
\begin{remark}
To prove an optimal error estimate for the stress in the $L^2$ norm, we  can follow the idea from \cite{Stenberg-1}  to use
 a mesh dependent norm technique.  In particular, this will lead to
 $$
 \|\sigma-\sigma_h\|_{L^2(\Omega)}\leq Ch^{k+1}|\sigma|_{H^{k+1}(\Omega)}.
 $$
\end{remark}

\section{First order mixed elements}
 In order to get  first order mixed elements, we propose to  take the following discrete displacement space
 \an{ \lab{Vh2}
   V_{1,h}: = \{v\in L^2(\Omega; \mathbb{R}^n),  \
         v|_K\in R(K)\ \hbox{ for any } K\in\mathcal{T}_h \}.
    }
 To  design the space for the stress, we define
 \begin{equation}
 \Sigma_{1,h}:=\{\tau\in H^1(\Omega;  \mathbb{S}), \tau|_K\in P_1(K, \mathbb{S})\text{ for any }K\in\mathcal{T}_h\}.
 \end{equation}
Since the pair $(\Sigma_{1,h}, V_{1,h})$ is unstable, we propose to enrich $\Sigma_{1,h}$ by some $n-1$ dimensional simplex bubble  function spaces.
Given a $n-1$ dimensional simplex $F$ of $\mathcal{T}_h$, let $\omega_F:=K^-\cup K^+$ denote the union of two elements that share $F$.  Define
\begin{equation}
\begin{split}
 \hat{\mathbb{B}}_F:=&\bigg\{\tau\in \hat{\Sigma}_{2,h}^\ast, \tau=0 \text{ on }\Omega\backslash\omega_F, \int_F \tau\nu\cdot p ds=0 \text{ for any } p\in \big(R(\omega_F)|_F\big)^\perp,\\
 & \text{ the values of } \int_K\tau:\theta d\b x \text{ vanish for any }\theta\in M_2(K), K=K^- \text{ and } K^+ \bigg\}.
\end{split}
\end{equation}
This allows for defining the following enriched stress space
 \begin{equation}
   \hat{\Sigma}_{1, h}^+=\Sigma_{1, h}+\sum\limits_{F}\hat{\mathbb{B}}_F.
 \end{equation}
 For this enriched space $\hat{\Sigma}_{1, h}^+$,   the number of degrees of freedom on each simplex is 18 and 48 for $n=2, 3$, respectively, which are the simplest conforming mixed elements so far. A  similar argument of Lemma \ref{lemma4.1}  shows that $\hat{\Sigma}_{1, h}^+$
 is a direct sum of $\Sigma_{1, h}$ and
$\sum\limits_{F}\hat{\mathbb{B}}_F$.

The mixed finite element approximation of Problem (1.1) reads: Find
   $(\sigma_h,~u_h)\in\hat{\Sigma}_{1,h}^+\times V_{1,h}$ such that
 \e{ \left\{ \ad{
    (A\sigma_h, \tau)+({\rm div}\tau, u_h)&= 0 &&
              \hbox{for all \ } \tau \in\hat{\Sigma}_{1,h}^+,\\
     (\d\sigma_h, v)& = (f, v) &&  \hbox{for all \ } v\in V_{1,h}.
      } \right. \lab{DP2}
    }
It follows from $\d  \hat{\Sigma}_{1,h}^+ \subset V_{1,h}$ that $\d  \tau=0$ for
   any $\tau\in Z_h$, which  implies the above K-ellipticity condition
	\meqref{below}.  A similar proof of Lemma \ref{lemma1}  shows the discrete inf--Sup  condition \ref{inf-sup}.
In particular, there exists an interpolation operator $I_h: H^1(\Omega, \mathbb{S})\rightarrow \hat{\Sigma}_{1,h}^+$ such that
\begin{equation}
\|\tau-I_h\tau\|_0+h\|\d(\tau-I_h\tau)\|\leq h^k \|\tau\|_k, k=1,2,
\end{equation}
and
\begin{equation}
\int_K \d (\tau-I_h\tau):p d\b x=\int_{\partial K}(\tau-I_h\tau)\nu\cdot p ds=0\text{ for any }p\in R(K)
\end{equation}
for any $K\in\mathcal{T}_h$.  A summary of these results leads  to the error  estimates in the following theorem.

\begin{theorem}\label{MainError2} Let
  $(\sigma, u)\in\Sigma\times V$ be the exact solution of
   problem \meqref{eqn1} and $(\tau_h, u_h)\in\hat{\Sigma}_{1,h}^+\times
   V_{1,h}$ the finite element solution of \meqref{DP2}.  Then,
\an{ \lab{t2} \|\sigma-\sigma_h\|_{H({\rm div})}
    + \|u-u_h\|_{0}&\le     Ch(\|\d \sigma\|_{1}+\|u\|_{1}),
      }
and
\begin{equation}
\|\sigma-\sigma_h\|_0\leq Ch^2\|\sigma\|_2.
\end{equation}
\end{theorem}

\section{The added face bubble  functions in three dimensions}
Let  $F:=\triangle_2\b x_1\b x_2 \b x_3$ be a face  of element $K:=\triangle_3\b x_0\b x_1 \b x_2\b x_3$,  we  construct
 the added face bubble functions. We have three face bubble functions of the Lagrange element of order $4$:
 \begin{equation}
 \varphi_{i, F}=\lambda_1\lambda_2\lambda_3(\lambda_i-\frac{1}{4}), i=1, 2, 3.
 \end{equation}
 Note that $\varphi_{i,F}$ vanish on face $F^\prime$ other than $F$ of $K$.

 Let $\b t_{i, F}$, $i=1, 2, 3$,  be unit tangential vectors of three edges of $F$. Let
 \begin{equation}
 T_{i, F}=\b t_{i, F}\b t_{i, F}^T, i=1, 2, 3.
 \end{equation}
 Define $T^\perp_{j, F}, j=1, 2, 3$  such that
 \begin{equation}
 T_{i,F}:T_{j, F}^\perp=0,  T^\perp_{j, F}:T^\perp_{l, F}=\delta_{jl},  i, j, l=1, 2, 3.
 \end{equation}
 This allows the definition of the following space
 \begin{equation}
\Sigma_{F, b}:=\text{span}\{\varphi_{i, F}T_{j, F}^\perp, i, j=1, 2, 3\}.
 \end{equation}
 On the face $F$, we have
 \begin{equation}
  x_i=\sum\limits_{j=1}^2(x_{i,j}-x_{i,3})\lambda_j+x_{i,3}, i=1,2,3,
 \end{equation}
 where $\b x_i=(x_{i,0}, x_{i,1}, x_{i,2}, x_{i,3})$, $i=0, 1, 2, 3$.  We need a basis of the restriction
 of the rigid motion space on the face $F$:
 \begin{equation}
 v_{1,F}=\p{1\\0 \\0}, v_{2,F}=\p{0\\1 \\0}, v_{3,F}=\p{0\\0 \\1},
 \end{equation}
and
\begin{equation}
 v_{4,F}=\p{(x_2-x_{2,F})\\-(x_1-x_{1,F}) \\0}, v_{5,F}=\p{0\\(x_3-x_{3,F}) \\-(x_2-x_{2,F})},
 v_{6,F}=\p{(x_3-x_{3,F})\\0 \\-(x_1-x_{1,F})}.
 \end{equation}
 Here $\b x_F=(x_{1, F}, x_{2, F}, x_{3, F})$ is the center of $F$. Define the basis $v_{i,F}^\perp$, $i=1, 2, 3$,  of the orthogonal complement space of the restriction  of the rigid motion space on the face $F$ with respect to $P_1(F, \mathbb{R}^3)$, such that
 \begin{equation}
 \int_F v_{i,F}^\perp\cdot v_{j, F}ds=0, i=1,2,3, j=1,\cdots, 6.
 \end{equation}
Then we define $\tau_{i, F}^\ast\in\Sigma_{F, b}$, $i=1, \cdots, 6$ such that
\begin{equation}
\frac{1}{|F|}\int_F \tau_{i, F}^\ast\nu_F\cdot v_{j,F}ds=\delta_{i,j}, j=1, \cdots, 6,
\text{ and }\int_F\tau_{i, F}^\ast\nu_F\cdot v_{k,F}^\perp ds=0, k=1,2,3.
\end{equation}
Finally, we take $\delta_{i, F}\in \Sigma_{K, 4,b}$ such that $\d \tau_{i, F}=\d (\tau_{i,F}^\ast+\delta_{i,F})\in P_1(K, \mathbb{R}^3)$. Then
\begin{equation}
\mathbb{B}_F^1=\text{span}\{\tau_{i, F}, i=1, \cdots, 6\}\text{ and }\mathbb{B}_F^2=\text{span}\{\tau_{i, F}, i=4, 5, 6\}.
\end{equation}
{\bf Example} Let $F=\triangle \b x_1\b x_2\b x_3$ with $\b x_1=(0,0,0)^T$, $\b x_2=(1,0,0)^T$, and $\b x_3=(0,1,0)^T$
 and $\nu_F=(0,0,1)^T$.  We have
 \begin{equation}
 T_{1,F}=\p{1&0&0\\ 0&0&0\\0&0&0 }, T_{2,F}=\p{1&-1&0\\ -1&1&0\\0&0&0 }, T_{3,F}=\p{0&0&0\\ 0&1&0\\0&0&0 }.
 \end{equation}
 This implies that
 \begin{equation}
 T_{1,F}^\perp=\p{0&0&0\\ 0&0&0\\0&0&1 }, T_{2,F}^\perp=\p{0&&1\\ 0&0&0\\1&0&0 }, T_{3,F}^\perp=\p{0&0&0\\ 0&0&1\\0&1&0 }.
 \end{equation}
 In addition,  a basis of the restriction
 of the rigid motion space on the face $F$ reads
 \begin{equation}
 v_{1,F}=\p{1\\0 \\0}, v_{2,F}=\p{0\\1 \\0}, v_{3,F}=\p{0\\0 \\1},
 \end{equation}
and
\begin{equation}
 v_{4,F}=\p{(x_2-\frac{1}{3})\\-(x_1-\frac{1}{3}) \\0}, v_{5,F}=\p{0\\0 \\-(x_2-\frac{1}{3})},
 v_{6,F}=\p{0\\0 \\-(x_1-\frac{1}{3})}.
 \end{equation}
Hence
\begin{equation}
 v_{1,F}^\perp=\p{(x_1-\frac{1}{3})\\0 \\0}, v_{2,F}^\perp=\p{0\\(x_2-\frac{1}{3})\\0 },
 v_{3,F}^\perp=\p{(x_2-\frac{1}{3})\\(x_1-\frac{1}{3}) \\0}.
 \end{equation}

\section{Numerical test}

We compute a 2D pure displacement problem  on the unit square $\Omega=[0,1]^2$
     with a homogeneous boundary condition
         that $u\equiv 0$ on $\partial\Omega$.
In the computation, we let the compliance tensor in \meqref{eqn1}
   \a{
      A \sigma &= \frac 1{2\mu} \left(
       \sigma - \frac{\lambda}{2\mu + n \lambda} \operatorname{tr}(\sigma)
        \delta \right), \quad n=2,  }
  where $\delta=\p{1 &0\\0&1}$, and $\mu=1/2$ and $\lambda=1$ are the
    Lam\'e constants.
Let  the  exact solution be
   \e{\lab{e-2}   u= \p{ e^{x-y}x (1-x) y (1-y)\\
    \sin(\pi x)\sin(\pi y)}. }
The true stress function $\sigma$
     and the load function $f$ are defined by the equations in
    \meqref{eqn1},   for the given  solution $u$.

In the computation, the level one grid consists of two right triangles,
   obtained by cutting the unit
   square with a north-east line.
Each grid is refined into a half-sized grid uniformly, to get
   a higher level grid.
In all the computation, the discrete systems of equations are
solved by Matlab backslash solver.

We use the bubble enriched $P_2$ symmetric stress finite element with
   $P_1$ discontinuous displacement finite element, $k=2$ in \meqref{Vh} and in
    \meqref{Sh-plus}, and $k=2$ in \meqref{Sh}.
That is, 3 $P_3$ bubbles are enriched each edge.
In Table \mref{b-2}, the errors and the convergence order
   in various norms are listed  for the true solution \meqref{e-2}.
The optimal order of convergence is observed  for both displacement and stress,
  see Table \mref{b-2}, as shown in the theorem.

\begin{table}[htb]
  \caption{  The errors, $\epsilon_h =\sigma -\sigma_h$,
     and the order of convergence, by the 2D $k=2$ element in
      \meqref{Sh-plus} and \meqref{Vh}, for \meqref{e-2}.}
\lab{b-2}
\begin{center}  \begin{tabular}{c|cc|cc|cc}  
\hline & $ \|u-u_h\|_{0}$ &$h^n$ &
    $ \|\epsilon_h\|_0$ & $h^n$  &
    \multispan{1} $ \|\d\epsilon_h \|_0$  & $h^n$ \\ \hline
1&  0.27452&0.0&  1.24637&0.0&   6.97007772&0.0 \\
2&  0.07432&1.9&  0.18054&2.8&   2.13781130&1.7 \\
3&  0.01959&1.9&  0.02429&2.9&   0.57734125&1.9 \\
4&  0.00497&2.0&  0.00314&2.9&   0.14709450&2.0 \\
5&  0.00125&2.0&  0.00040&3.0&   0.03694721&2.0 \\
      \hline
\end{tabular}\end{center} \end{table}

As a comparison,  we also test   the Arnold--Winther element from \cite{Arnold-Winther-conforming}, which  has a same degree of freedom as ours, 21, on each  element.
But the displacement in that element
   is approximated by the rigid-motion space only,
   instead of the full $P_1$ space, i.e., 3 dof vs 6 dof on each triangle. The total degrees of freedom for the stress for the new element are $3|\mathbb{V}|+3|\mathbb{E}|+3|\mathbb{K}|$, where    $|\mathbb{V}|$, $|\mathbb{E}|$, and $|\mathbb{K}|$ are the numbers of vertices, edges and elements of $\mathcal{T}_h$, respectively, while those for the Arnold--Winther element are $3|\mathbb{V}|+4|\mathbb{E}|$.  Since the
       three bubble functions on each element can be easily condensed,  these two elements almost have  the same complexity for solving.
The errors and the orders of convergence are listed in Table \mref{b-3}.
Because the new element uses the full $P_1$ displacement space, the order of
convergence is one higher than that of the Arnold--Winther element.
Also as the new element includes the full $P_2$ stress space, the order of
  convergence of stress is one order higher, see  the data in Tables \mref{b-2} and \mref{b-3}.

\begin{table}[htb]
  \caption{  The errors, $\epsilon_h =\sigma -\sigma_h$,
     and the order of convergence, by the
     Arnold-Winther 21/3 element\cite{Arnold-Winther-conforming}, for \meqref{e-2}.}
\lab{b-3}
\begin{center}  \begin{tabular}{c|cc|cc|cc}  
\hline & $ \|u-u_h\|_{0}$ &$h^n$ &
    $ \|\epsilon_h\|_0$ & $h^n$  &
    \multispan{1} $ \|\d\epsilon_h \|_0$  & $h^n$ \\ \hline
 1&  0.30554&0.0&  1.58058&0.0&  10.31991249&0.0 \\
 2&  0.22589&0.4&  0.89927&0.8&   6.81340378&0.6 \\
 3&  0.10922&1.0&  0.25584&1.8&   3.61633797&0.9 \\
 4&  0.05354&1.0&  0.06633&1.9&   1.83690959&1.0 \\
 5&  0.02661&1.0&  0.01674&2.0&   0.92212628&1.0 \\      \hline
\end{tabular}\end{center} \end{table}

\end{document}